\documentclass{amsart}
\usepackage{amssymb,latexsym}
\usepackage{graphicx}
\usepackage{color}
\usepackage{enumitem}
\usepackage{cancel}

\newtheorem{theorem}{Theorem}[section]

\newtheorem{remark}[theorem]{Remark}

\newtheorem{conjecture}[theorem]{Conjecture}
\newtheorem{corollary}[theorem]{Corollary}
\newtheorem{proposition}[theorem]{Proposition}

\def\daniele#1 {\fbox {\footnote {\ }}\ \footnotetext { From Daniele: {\color{blue}#1}}}
\def\massimo#1 {\fbox {\footnote {\ }}\ \footnotetext { From Massimo: {\color{red}#1}}}

\begin{document}

\title{Permutation polynomials, fractional polynomials, and algebraic curves}
\author{Daniele Bartoli\textsuperscript{\,1,$\dagger$}}
\author{Massimo Giulietti\textsuperscript{\,2}}
\address{\textsuperscript{1}Department of Mathematics and Computer Science, University of Perugia, 06123 Perugia, Italy}
\email{daniele.bartoli@unipg.it}
\address{\textsuperscript{2}Department of Mathematics and Computer Science, University of Perugia, 06123 Perugia, Italy}
\email{massimo.giulietti@unipg.it}
\address{\textsuperscript{$\dagger$}Corresponding Author}

\date{\today}

\maketitle

\begin{abstract}
In this note we prove a conjecture by Li, Qu, Li, and Fu 
on  permutation trinomials over $\mathbb{F}_3^{2k}$. 
In addition, new examples and generalizations of some families of permutation polynomials of $\mathbb{F}_{3^k}$ and $\mathbb{F}_{5^k}$ are  given. We also study permutation quadrinomials of type $Ax^{q(q-1)+1} + Bx^{2(q-1)+1} + Cx^{q} + x$. 
Our method is based on the investigation of an algebraic curve associated with a {fractional polynomial} over a finite field.  
\end{abstract}

\section{Introduction}

Let $q=p^h$ be a prime power. A polynomial $f(x)\in \mathbb{F}_q[x]$ is a {\it permutation polynomial} (PP) if it is a bijection of the finite field $\mathbb{F}_q$ into itself. On the other hand, each permutation of $\mathbb{F}_q$ can be expressed as a polynomial over $\mathbb{F}_q$. Permutation polynomials were first studied by Hermite and  Dickson; see \cite{Dickson1896,Hermite1863}.

In general it is not difficult to construct a random PP for a given field $\mathbb{F}_q$. Particular, simple structures or additional extraordinary properties are usually required by applications of PPs in other areas of mathematics and engineering, such as cryptography, coding theory, or combinatorial designs. Permutation polynomials meeting these criteria are usually difficult to find. For a deeper introduction on the connections of PPs with other fields of mathematics we refer to \cite{MuPa, Hou2015} and the references therein. 

In this work we deal with a particular class of PP. For a prime $p$ and a positive integer $m$, let  $\mathbb{F}_{p^m}$ be the finite field with $p^m$ elements. Given a polynomial $h(x) $ over  $\mathbb{F}_{p^m}$, a divisor $d$ of $p^m-1$, and an integer $r$ with  $1 \leq r<(p^m-1)/d$, let
$$f_{r,d,h}(x)=x^rh\left( x^{\frac{p^m-1}{d}}\right).$$
A useful criterion to decide weather $f_{r,d,h}$ permutes $\mathbb{F}_{p^m}$ is the following.

\begin{theorem}\label{ThPLZ} \cite{PL2001,Zieve2009}
The polynomial $f_{r,d,h}(x)$ is a PP of $\mathbb{F}_{p^m}$ if and only if $gcd(r, (p^m-1)/d) =1$ and $x^rh(x)^{(p^m-1)/d}$ permutes the set $\mu_d$ of the $d$-th roots of unity in $\mathbb{F}_{p^m}$.
\end{theorem}

{Let $q=p^n$. For $h(x) =\sum_{i=0}^{\ell} a_i x^i$ a polynomial over $\mathbb F_{q^2}$, by Theorem \ref{ThPLZ} $x^rh\left( x^{q-1}\right)$ permutes $\mathbb{F}_{q^2}$ if and only if $x^rh\left( x\right)^{q-1}$ permutes $\mu_{q+1}$. If this is the case and in addition $h(x) \in \mathbb{F}_{q}[x]$, then for each  $z\in \mu_{q+1}$ we have that  
$$z^rh\left( z\right)^{q-1}=z^r\frac{\left(h( z)\right)^q}{h(z)}=z^r\frac{h(1/z)}{h(z)}=z^{r-\ell}\frac{\widetilde{h}(z)}{h(z)},$$
 where $\deg(h)=\ell$ and   $\widetilde{h}(x) =\sum_{i=0}^{\ell} a_{\ell-i} x^i$. We call the rational function $x^{r-\ell}\frac{\widetilde{h}(x)}{h(x)}$ the \emph{fractional polynomial} associated with the PP $x^rh\left( x^{q-1}\right)$. Conversely, given a fractional polynomial $x^{r-\ell}\frac{\widetilde{h}(x)}{h(x)}$ which permutes $\mu_{q+1}$ we call $x^rh\left( x^{q-1}\right)$ the associated permutation polynomial.}  

A standard approach to the problem of deciding whether a polynomial $f(x)$ is a PP 
is the investigation of the plane algebraic curve 
$$
\mathcal C_f:\frac{f(x)-f(y)}{x-y}=0;
$$
in fact, $f$ is a PP over $\mathbb F_{p^m}$ if and only if $\mathcal C_f$ has no
 $\mathbb F_{p^m}$-rational point $(a,b)$ with $a\neq b$. In the case where $p^m=q^2$ and
  $f$ is of type $f_{r,q+1,h}$ with $h\in \mathbb F_q[x]$, it can be more effective to study
   the curve, with degree lower than $\mathcal C_f$, defined by the equation
   
$$
x^{r-\ell}\frac{\widetilde{h}(x)}{h(x)}-y^{r-\ell}\frac{\widetilde{h}(y)}{h(y)}=0
$$
and check whether it has some $\mathbb F_{q^2}$-rational points $(a,b)$ with 
$a\neq b$ and $a^{q+1}=b^{q+1}=1$. 

In this note we use both such methods to solve some conjectures on permutation polynomials, shorten several proofs that have recently appeared in the literature, and provide some new examples of permutation quadrinomials and permutation functions of type $\alpha(x)/\beta(x)$, where $\alpha,\beta \in \mathbb{F}_{q}[x]$.

\section{On some conjectures on permutation polynomials}\label{Sec2}

In \cite{LQLF2016} the authors presented the following conjecture about permutation trinomials in characteristic $3$.

\begin{conjecture}\cite[Conjecture 5.1]{LQLF2016}
\begin{enumerate}
\item Let $q=3^k$, $k$ even, and $$f(x) =x^{\ell q+\ell+5}+x^{(\ell+5)q+\ell}-x^{(\ell - 1)q+\ell+6},$$ where $gcd(5 +2\ell, q-1) =1$. Then $f(x)$ is a permutation trinomial over $\mathbb{F}_{q^2}$.
\item Let $q=3^k$, $$f(x) =x^{\ell q+\ell+1}-x^{(\ell+4)q+\ell-3}-x^{(\ell - 2)q+\ell+3}$$ and $gcd(1 +2\ell , q-1) =1$. Then $f(x)$ is a permutation trinomial over $\mathbb{F}_{q^2}$.
\item Let $q=3^k$, $$f(x) =x^{\ell q+\ell+1}+x^{(\ell+2)q+\ell-1}-x^{(\ell - 2)q+\ell+3}$$ and $gcd(1 +2\ell, q-1) =1$. Then $f(x)$ is a permutation trinomial over $\mathbb{F}_{q^2}$ if $k\not \equiv 2 \pmod{4}$.
\end{enumerate}

\end{conjecture}

It has been noticed in \cite{Li2017} that by Theorem \ref{ThPLZ} the above conjecture can be rephrased as follows.

\begin{conjecture}\cite[Conjecture 5.2]{LQLF2016} \label{Conj}

\begin{enumerate}
\item Let $q=3^k$, $k$ even and $g(x) =\frac{-x^7+x^6+x}{x^6+x-1}$. Then $g(x)$ permutes $\mu_{q+1}$.
\item Let $q=3^k$ and $g(x) =\frac{x^6+x^4-1}{-x^7+x^3+x}$. Then $g(x)$ permutes $\mu_{q+1}$.
\item Let $q=3^k$ and $g(x) =\frac{-x^5+x^3+x}{x^4+x^2-1}$. Then $g(x)$ permutes $\mu _{q+1}$ if $k \not \equiv 2\pmod{4}$.
\end{enumerate}
\end{conjecture}

In \cite{Li2017} the author settles Conjecture \ref{Conj}(2) and Conjecture \ref{Conj}(3) by determining some quadratic factors of a five-degree polynomial and a seven-degree polynomial. 

The aim of this section is twofold. On yhe one hand,  we settle Conjecture \ref{Conj}(1) using different arguments with the respect to those of  \cite{Li2017}. On the other hand, we show that basic tools from Algebraic Geometry can be very useful when dealing with permutation polynomials (see also \cite{BGZ2016,BGQZ2017}); in particular they can provide shorter and less technical proofs.


\begin{theorem}
Let $q=3^k$, $k$ even. The function $g(x)=\frac{-x^7+x^6+x}{x^6+x-1}$ permutes $\mu_{q+1}$.
\end{theorem}

\proof

Note that if $x \in \mu_{q+1}$, then 
$$\left(\frac{-x^7+x^6+x}{x^6+x-1}\right)^{q}=\left(-x\frac{x^6-x^5-1}{x^6+x-1}\right)^{q}=\frac{-1}{x}\left(\frac{-x^6-x+1}{-x^6+x^5+1}\right)=\left(\frac{x^6+x-1}{-x^7+x^6+x}\right),$$
that is $g(\mu_{q+1}) \subset \mu_{q+1}$. We have only to show that $g(x)$ is injective on $\mu_{q+1}$. 

Let $\mathcal{C}$ be the plane curve given by the affine equation $H(x,y)=0$, where 
$$H(x,y)=(-x^7+x^6+x)(y^6+y-1) -(x^6+x-1)(-y^7+y^6+y)=0.$$

The function $g(x)$ permutes  $\mu_{q+1}$ if and only if the curve $\mathcal{C}$ does not have  any point $(\overline x,\overline y)\in \mu_{q+1}^2$ off the line $x=y$.

The polynomial $H$ factorizes as 
$$H(x,y)=(x-y) F(x,y)G(x,y),$$
where $$
\begin{array}{ll}
F(x,y)=& x^3y^3 + \omega ^2x^3y^2 + \omega ^5x^3y + \omega ^5x^3 + \omega ^6x^2y^3+ x^2y^2+ 2x^2y + \omega ^5x^2\\ 
 &        + \omega ^7xy^3 + 2xy^2 + xy + \omega ^2x + 
        \omega ^7y^3 + \omega ^7y^2 + \omega ^6y + 1\\
\end{array}$$
and 
$$\begin{array}{ll}
G(x,y)=& x^3y^3 + \omega ^6x^3y^2 + \omega ^7x^3y + \omega ^7x^3 + \omega ^2x^2y^3 + x^2y^2
        + 2x^2y+ \omega ^7x^2  \\
        &+ \omega ^5xy^3 + 2xy^2 + xy + \omega ^6x+ 
        \omega ^5y^3 + \omega ^5y^2 + \omega ^2y + 1,\\
        \end{array}$$
for some primitive element $\omega$ of $\mathbb F_{9}$. Note that both $F(x,y)$ and $G(x,y)$ are polynomials defined over $\mathbb{F}_{3^k}$, since $k$ is even.

Suppose that $\overline x, \overline y \in \mu_{q+1}$ are such that $F(\overline x,\overline y)=0$. Then $0=\overline{x}^3\overline{y}^3F(\overline x,\overline y)^{q}=\overline{x}^3\overline{y}^3 F(1/\overline x,1/\overline y)=$

$$
\begin{array}{l}
\overline x^3\overline y^3 +  \omega^6\overline x^3\overline y^2 +  \omega^7\overline x^3\overline y +  \omega^7\overline x^3 +  \omega^2\overline x^2 \overline y^3 +\overline x^2\overline y^2 + 
    2 \overline x^2\overline y +  \omega^7\overline x^2\\ +  \omega^5\overline x\overline y^3 + 2\overline x\overline y^2 + \overline x \overline y  +  \omega^6\overline  x +  \omega^5\overline y^3 + 
     \omega^5\overline y^2 +  \omega^2\overline y + 1.
     \end{array}$$

The resultant between $F(x,y)$ and $x^3y^3F(1/x,1/y)$ with the respect to $y$ is $(x + \omega^2)^9(x + \omega^6)^9$. This implies that the common points of the curves with equations $F(x,y)=0$ and $x^3y^3F(1/x,1/y)=0$ belong to the lines $x=2\omega^2,2\omega^6$, but in this case $x^q=x\neq 1/x$. Therefore no points $(\overline x,\overline y) \in \mu_{q+1}^2$ satisfy $F(\overline x,\overline y)=0$. A similar argument holds for $G(x,y)$. Then the function $g(x)=\frac{-x^7+x^6+x}{x^6+x-1}$ permutes $\mu_{q+1}$.

\endproof

\begin{remark}
Similar arguments can be used to settle the Conjectures {\rm \ref{Conj}(2)} and {\rm \ref{Conj}(3)}. In the former case, the curve splits into six conic components, apart from $x-y=0$, defined by  
$$x y +  \eta^2 x +  \eta^{15} y + 1=0, \qquad x y +  \eta^5 x +  \eta^{18} y + 1=0, \qquad x y +  \eta^6 x +  \eta^{19} y + 1=0$$
$$x y +  \eta^{15} x +  \eta^2 y + 1=0, \qquad x y +  \eta^{18} x +  \eta^5 y + 1=0, \qquad x y +  \eta^{19} x +  \eta^6 y + 1=0,$$
where $\eta$ is some primitive element of $\mathbb{F}_{27}$. 

In the latter case the curve splits into four conic components defined over $\mathbb{F}_{81}$, namely 
$$x y +  \xi^{25} x +  \xi^{65} y + 1=0,\qquad x y +  \xi^{35} x +  \xi^{75} y + 1=0,$$
$$x y +  \xi^{65} x +  \xi^{25} y + 1=0, \qquad x y +  \xi^{75} x +  \xi^{35} y + 1=0;$$
here $\xi$ is a primitive element of $\mathbb F_{81}$.

In both cases it can be easily proved that no points $(x,y) \in \mu_{q+1}^2$ with $x\neq y$ belong to a component of the curve, and therefore the corresponding fractional functions $g(x)$ permute $\mu_{q+1}$. 
\end{remark}

We now provide alternative and shorter proofs of two results from \cite{MG2017}, referring to two conjectures presented in \cite{WL2017}.
\begin{proposition}\label{prop25}
Let $q = 5^k$, $k$ even. Then $g(x) = -x\frac{(x^2-2)^2}{(x^2+2)^2}$ permutes $\mu_ {q+1}$.
\end{proposition}
\proof
It is enough to prove that $g(x)$ is injective over $\mu_{q+1}$. The corresponding plane algebrai  curve is
$$\mathcal{C} :  \frac{-x(x^2-2)^2(y^2+2)^2 +y(y^2-2)^2(x^2+2)^2}{x-y} =0.$$
The equation of $\mathcal C$ can be rewritten as 
$$
F_1(x,y)\cdot F_2(x,y)=0
$$
where
$$
F_1(x,y)=x^2 y^2 +  \omega^3 x^2 y + 3 x^2 +  \omega^3 x y^2 + x y +  \omega^{15} x + 3 y^2 + 
         \omega^{15} y + 1
$$
and
$$         
    F_2(x,y)=     
     x^2 y^2 +  \omega^{15} x^2 y + 3 x^2 +  \omega^{15} x y^2 + x y +  \omega^3 x + 3 y^2 + 
         \omega^3 y + 1,$$
for a primitive element $\omega$ of $\mathbb{F}_{25}$. A point $(\overline x, \overline y) \in \mu_{q+1}^2$ belongs to a component $F_i(x,y)=0$ if and only if  satisfies both $F_i(\overline x,\overline y)=F_i(1/\overline x,1/\overline y)=0$. Now the resultant of the two polynomials $F_i(x,y)$ and $x^2y^2F_i(1/x,1/y)$ with the respect to $y$ is always
$R(x)=4 x^8 + 4 x^6 + 3 x^4 + 4 x^2 + 4$. The polynomial $R(x)$ factorizes as $\prod _{j \in J} (x+\omega^j)$, where $J=\{3,4,8,9,15,16,20,21\}$. None of its roots is in $\mu_{q+1}$ and therefore  $\mathcal{C} $ does not contains points of type $(\overline x, \overline y) \in \mu_{q+1}^2$; that is, $g(x)$ permutes $\mu_{q+1}$.
\endproof

\begin{proposition}\label{Prop:F5}
Let $q = 5^k$, $k$ odd. Then $g(x) = x\frac{(x^2-x+2)^2}{(x^2+x+2)^2}$ permutes $\mathbb{F}_q$.
\end{proposition}
\proof
It is enough to show that $g(x)$ is injective over $\mathbb{F}_q$. The corresponding curve is
$$\mathcal{C} \ : \ \frac{x(x^2-x+2)^2(y^2+y+2)^2 -y(y^2-y+2)^2(x^2+x+2)^2}{x-y} =0.$$
This curve as an (affine) $\mathbb{F}_q$-rational point $(\overline x,\overline y)$, $\overline x\neq \overline y$, if and only if $g(x)$ is a permutation of $\mathbb{F}_q$.
The equation of $\mathcal C$ can be written as $F_1(x,y)F_2(x,y)=0$, where
$$F_1(x,y)= (x^2 y^2 +  \omega^7 x^2 y +  \omega^2 x^2 +  \omega^7 x y^2 +  \omega^3 x y +  \omega x + 
         \omega^2 y^2 +  \omega y + 4)$$
         
$$F_2(x,y)=(x^2 y^2 +  \omega^{11} x^2 y +  \omega^{10} x^2 +  \omega^{11} x y^2 +  \omega^{15} x y +  \omega^5 x +
         \omega^{10} y^2 +  \omega^5 y + 4),$$
for some primitive element $\omega $  of  $\mathbb{F}_{25}$. Note that the curves $F_i(x,y)=0$ are not defined over $\mathbb{F}_q$, since $k$ is odd. The unique $\mathbb{F}_{q}$-rational points of $\mathcal{C}$ satisfy $F_{1}(x,y)=F_2(x,y)=0$.

The resultant between $F_1(x,y)$ and $F_2(x,y)$ with the respect to $y$ is 
$$R(x)= 4 x^8 + x^7 + 4 x^6 + x^5 + 3 x^4 + 3 x^3 + x^2 + 2 x + 4$$

and its roots are $\omega^i$ , $i \in \{1,4,5,13,14,17,20,22\}$. None of them is in $\mathbb{F}_{5^k}$ since $k$ is odd and the claim follows.

\endproof


Using the same argument one can easily prove \cite[Theorem 3.2]{LQW2017} and \cite[Theorem 3.25.]{LQW2017}. Namely,
when $q=2^k$ the fractional polynomials 
$$\frac{x^8+x^7+x^6+x^5+x^3+x^2+x}{x^7+x^6+x^5+x^3+x^2+x+1}, \qquad \frac{x^8+x^7+x^5+x^3+x}{x^7+x^5+x^3+x+1}$$
permute $\mu_{q+1}$. Let $\mathbb{F}_{2^7}=\langle w \rangle$, where $w^7 + w + 1=0$.  The equations of the associated curves can be written as
$$
\prod_{i=0}^{6}(xy + w^{19\cdot 2^i} x + w^{19\cdot 2^i} y + 1)=0$$
and 
$$\prod_{i=0}^{6}(xy + w^{3\cdot 2^i} x + w^{3\cdot 2^i} y + 1)=0,$$
respectively. Each component can be written as $y=\frac{\xi +w}{\xi x+1}$, for some $\xi \in \mathbb{F}_{2^7}$. It is easily seen that if $q=2^n$ with $n\equiv 0 \pmod 7$ then $y^q=1/y$ for any $x\in \mu_{q+1}$, since $\xi\in \mathbb{F}_q$. On the other hand if $n\not\equiv 0 \pmod 7$ then $\xi\notin \mathbb{F}_{q}$ and  $y^q=1/y$ if only if $y=x=1$. Therefore both the fractional polynomials permute $\mu_{q+1}$ if and only if $q=2^n$ with $n\not\equiv 0 \pmod 7$.

We end this section by discussin briefly two conjectures about permutation trinomials presented by Gupta and Sharma  \cite{GS2016} and proved in  \cite{ZHF2017, WYDM2017}:

\begin{itemize}
\item the polynomial $f(x) =x^5+x^{3\cdot 2^m+2} +x^{4\cdot 2^m+1} \in \mathbb{F}_{2^{2m}}[x]$ is a permutation trinomial over $\mathbb{F}_{2^{2m}}$ if and only if $m \equiv 2\pmod 4$;
\item the polynomial $f(x) =x^5+x^{2^m+4} +x^{5\cdot 2^m} \in \mathbb{F}_{2^{2m}}[x]$ is a permutation trinomial over $\mathbb{F}_{2^{2m}}$ if and only if $m \equiv 2\pmod 4$.
\end{itemize}

 The proofs given in  \cite{ZHF2017, WYDM2017} rely on the fact that the rational functions 
$$\frac{x+x^2+x^5}{1+x^3+x^4},\qquad  \frac{1+x+x^3}{1+x^2+x^3},\qquad  \frac{1+x+x^5}{1+x^4+x^5}$$
or their inverses permute $\mu_{2^m+1}$. This follows quite easily from considering the associated curves. For instance, the curve associated with the function $\frac{1+x+x^5}{1+x^4+x^5}$ is 
 $$\frac{(x^5+x^4+x) (y^4+y+1) + (y^5+y^4+y) (x^4+x+1)}{x+y}=0,$$
 which splits as
 $$(x  y + \omega  x + \omega^4  y + 1)(x  y + \omega^2  x + \omega^8  y + 1)(x  y + \omega^4  x + \omega  y + 1)(x  y + \omega^8  x + \omega^2  y + 1)=0,$$
for some primitive element $\omega $ of  $\mathbb{F}_{16}$; this means that the components of the curve are not defined over $\mathbb{F}_{2^m}$ if $m \equiv 2 \pmod 4$. Whence, it is easy to show that no pairs $(x,y) \in \mu_{2^m+1}$, $x \neq y$ belong to such a curve.

\section{New examples of fractional permutations}

In this section we extend some results of Section \ref{Sec2} to  other fractional polynomials.
We begin with the case $p=5$.

\begin{proposition}
Let $q = 5^k$, $k$ odd, and $A,B,C,D\in \mathbb{F}_5$.  Then $g_{A,B,C,D}(x) = x\frac{(x^2+Ax+B)^2}{(x^2+Cx+D)^2}$ permutes $\mathbb{F}_q$ in all the cases listed in Table \ref{Table}.
\end{proposition}

\proof
In all the cases listed in Table \ref{Table} $\omega$ stands for a primitive element of $\mathbb F_{25}$.  The curve corresponding to the fractional polynomial has two components $F_1(x,y)=0$ and $F_2(x,y)=0$ not defined over $\mathbb{F}_{q^k}$, but over $\mathbb{F}_{q^{2k}}$. One can argue as in the proof of Proposition \ref{prop25}. The last column of the table indicates the roots of the resultant $R(x)$ of $F_1(x,y)$ and $F_2(x,y)$ with the respect to $y$. It is easily seen that in all these cases such roots do not belong to $\mathbb{F}_{q^k}$ since $k$ is odd; therefore,  the curve associated with $g_{A,B,C,D}$ does not have any $\mathbb{F}_{q^k}$-rational points with distinct coordinates. 

\endproof

\begin{footnotesize}

\begin{table}
\caption{Permutation polynomials of type $x\frac{(x^2+Ax+B)^2}{(x^2+Cx+D)^2}$ of $\mathbb{F}_{5^k}$, $k$ odd, $A,B,C,D \in \mathbb{F}_5$}\label{Table}
\tabcolsep = 0.1 mm
\begin{tabular}{|c|c|c|}
\hline
$A,B,C,D$ & $F_1(x,y)$ and $F_2(x,y)$ & Roots $\omega^i$ of $R(x)$\\
\hline
\hline
 $[ 4, 2 , 1, 2 ]$
 &

\begin{tabular}{c}
$x^2 y^2 + \omega^7 x^2 y + \omega^2 x^2 + \omega^7 x y^2 +  \omega^3 x y +  \omega x +  \omega^2 y^2 
    +  \omega y + 4$,\\
$x^2 y^2 +  \omega^{11} x^2 y +  \omega^{10} x^2 +  \omega^{11} x y^2 +  \omega^{15} x y +  \omega^5 x + 
     \omega^{10} y^2 +  \omega^5 y + 4$\\
\end{tabular}
&
$i=1,2,5,8,10,13,16,17$\\
\hline

$
    [ 4, 2 , 3, 3 ]
    $
    &

\begin{tabular}{c}
$x^2 y^2 +  \omega x^2 y +  \omega^{14} x^2 +  \omega x y^2 +  \omega^{16} x y +  \omega^{13} x + 
     \omega^{14} y^2 +  \omega^{13} y + 1$,\\
$x^2 y^2 +  \omega^5 x^2 y +  \omega^{22} x^2 +  \omega^5 x y^2 +  \omega^8 x y +  \omega^{17} x + 
     \omega^{22} y^2 +  \omega^{17} y + 1$\\
     \end{tabular}
&
$i=1,4,5,7,8,11,16,20$\\
\hline

$
    [ 2, 3 , 4, 2 ]
$
&
\begin{tabular}{c}
$x^2 y^2 +  \omega^{19} x^2 y +  \omega^2 x^2 +  \omega^{19} x y^2 +  \omega^4 x y +  \omega^{19} x + 
     \omega^2 y^2 +  \omega^{19} y + 1$,\\
$x^2 y^2 +  \omega^{23} x^2 y +  \omega^{10} x^2 +  \omega^{23} x y^2 +  \omega^{20} x y +  \omega^{23} x + 
     \omega^{10} y^2 +  \omega^{23} y + 1$\\
     \end{tabular}
&
$i=1,2,5,10,14,19,22,23$\\
\hline

$
    [ 2, 3 , 3, 3 ]
$
&
\begin{tabular}{c}
$x^2 y^2 +  \omega x^2 y +  \omega^{14} x^2 +  \omega x y^2 +  \omega^{15} x y +  \omega^7 x +  \omega^{14} y^2
    +  \omega^7 y + 4$,\\
$x^2 y^2 +  \omega^5 x^2 y +  \omega^{22} x^2 +  \omega^5 x y^2 +  \omega^3 x y +  \omega^{11} x + 
     \omega^{22} y^2 +  \omega^{11} y + 4$\\
     \end{tabular}
&
$i=2,4,7,10,11,19,20,23 $\\
\hline

$
    [ 1, 2 , 4, 2 ]
$
&
\begin{tabular}{c}
$x^2 y^2 +  \omega^{19} x^2 y +  \omega^2 x^2 +  \omega^{19} x y^2 +  \omega^3 x y +  \omega^{13} x + 
     \omega^2 y^2 +  \omega^{13} y + 4$,\\
$x^2 y^2 +  \omega^{23} x^2 y +  \omega^{10} x^2 +  \omega^{23} x y^2 +  \omega^{15} x y +  \omega^{17} x + 
     \omega^{10} y^2 +  \omega^{17} y + 4$\\
     \end{tabular}
&
$i=1,4,5,13,14,17,20,22 $
\\
\hline

$
    [ 1, 2 , 2, 3 ]$
&
\begin{tabular}{c}
$x^2 y^2 +  \omega^{13} x^2 y +  \omega^{14} x^2 +  \omega^{13} x y^2 +  \omega^{16} x y +  \omega x + 
     \omega^{14} y^2 +  \omega y + 1$,\\
$x^2 y^2 +  \omega^{17} x^2 y +  \omega^{22} x^2 +  \omega^{17} x y^2 +  \omega^8 x y +  \omega^5 x + 
     \omega^{22} y^2 +  \omega^5 y + 1$\\
     \end{tabular}
&
$i=4,8,13,16,17,19,20,23$
\\

\hline

$
    [ 3, 3 , 1, 2 ]
$
&
\begin{tabular}{c}
$x^2 y^2 +  \omega^7 x^2 y +  \omega^2 x^2 +  \omega^7 x y^2 +  \omega^4 x y +  \omega^7 x + 
     \omega^2 y^2 +  \omega^7 y + 1$,\\
$x^2 y^2 +  \omega^{11} x^2 y +  \omega^{10} x^2 +  \omega^{11} x y^2 +  \omega^{20} x y +  \omega^{11} x + 
     \omega^{10} y^2 +  \omega^{11} y + 1$\\
     \end{tabular}
&
$i=2,7,10,11,13,14,17,22$
\\
\hline

$
    [ 3, 3 , 2, 3 ]
$
&
\begin{tabular}{c}
$x^2 y^2 +  \omega^{13} x^2 y +  \omega^{14} x^2 +  \omega^{13} x y^2 +  \omega^{15} x y +  \omega^{19} x + 
     \omega^{14} y^2 +  \omega^{19} y + 4$,\\
$x^2 y^2 +  \omega^{17} x^2 y +  \omega^{22} x^2 +  \omega^{17} x y^2 +  \omega^3 x y +  \omega^{23} x + 
     \omega^{22} y^2 +  \omega^{23} y + 4$\\
     \end{tabular}
&
$i=7,8,11,14,16,19,22,23$
\\
\hline

$
    [ 0, 3 , 0, 2 ]$
&
\begin{tabular}{c}
$x^2 y^2 +  \omega^3 x^2 y + 3 x^2 +  \omega^3 x y^2 + x y +  \omega^{15} x + 3 y^2 +  \omega^{15} y 
    + 1$,\\
$x^2 y^2 +  \omega^{15} x^2 y + 3 x^2 +  \omega^{15} x y^2 + x y +  \omega^3 x + 3 y^2 +  \omega^3 y 
    + 1$\\
    \end{tabular}
&
$i=3,4,8,9,15,16,20,21$
\\

\hline

$[ 0, 2 , 0, 3 ]$&
\begin{tabular}{c}
$x^2 y^2 +  \omega^9 x^2 y + 2 x^2 +  \omega^9 x y^2 + 4 x y +  \omega^9 x + 2 y^2 +  \omega^9 y 
    + 1$,\\
$x^2 y^2 +  \omega^{21} x^2 y + 2 x^2 +  \omega^{21} x y^2 + 4 x y +  \omega^{21} x + 2 y^2 + 
     \omega^{21} y + 1$\\
   \end{tabular}
   &
   $i=2,3,9,10,14,15,21,22$
\\

\hline
\end{tabular}
\end{table}

\end{footnotesize}

Now we turn our attention to the case $\mathbb{F}_{3^k}$. We consider functions of type 
$$g_{A_1,A_2,A_3,A_4,A_5,A_6}(x)=\frac{-x^7+A_1x^6+A_2x^5+A_3 x^4+A_4x^3+A_5x^2+A_6x)}{(A_6x^6+A_5x^5+A_4x^4+A_3x^3+A_2x^2+A_1x-1)},$$
for $A_i \in \mathbb{F}_3$ and $A_1\neq 0$. In Table  \ref{Table2} 
$\omega$ stands for a primitive element of $\mathbb F_9$.
The components of the associated curve over $\mathbb{F}_{9}$ are indicated.  Note that every quadratic component  in the table is absolutely irreducible. It is easily seen that their only $\mathbb{F}_3$-rational points  lie also  on the line $x=y$. and therefore if $k$ is odd these conics do not have rational points off the line $x=y$. It is worth noting that all the conics for $k$ odd have points in $\mu_{3^k+1}^2$ and therefore the corresponding $g(x)$ does not permute $\mu_{3^k+1}$.

\begin{table}
\caption{Factorization over $\mathbb{F}_{3^2}$ of the curve associated with the function $g_{A_1,A_2,A_3,A_4,A_5,A_6}(x)$}\label{Table2}
\tabcolsep = 0.1 mm
\begin{tabular}{|c|c|c|c|}
\hline
$A_1,A_2,A_3,A_4,A_5,A_6$ & Factors &$A_1,A_2,A_3,A_4,A_5,A_6$ & Factors \\
\hline
\hline

$[ 1, 1, 0, 2, 1, 1 ]$&

$
\begin{array}{l}
      y + \omega, y + \omega^3 \\
      y + \omega^5, y + \omega^7 \\
      x + \omega, x + \omega^3 \\
      x + \omega^5, x + \omega^7 \\

      x y + \omega x + \omega^3 y + 1 \\
      x y + \omega^3 x + \omega y + 1 \\
\end{array}$&
$[ 1, 2, 2, 2, 0, 1 ]$&

$
\begin{array}{l}
      y^2 + \omega y + 1 \\
      y^2 + \omega^3 y + 1 \\

      x y + \omega^5 x + \omega^7 y + 1 \\
      x y + \omega^7 x + \omega^5 y + 1 \\
      x^2 + \omega x + 1 \\
      x^2 + \omega^3 x + 1 \\
\end{array}$

\\

\hline

$[ 2, 2, 1, 1, 2, 1 ]$&

$
\begin{array}{l}
      y + \omega^2,  y + \omega^6\\
      x + \omega^2, x + \omega^6\\

      x y + \omega^5 x + \omega^7 y + 1 \\
      x y + \omega^7 x + \omega^5 y + 1 \\
\end{array}$
&
$[ 1, 2, 2, 1, 1, 1 ]$&

$
\begin{array}{l}
      y + \omega^2,  y + \omega^6\\
      x + \omega^2,  x + \omega^6\\

      x y + \omega x + \omega^3 y + 1 \\
      x y + \omega^3 x + \omega y + 1 \\
\end{array}$
\\

\hline
$[ 2, 2, 1, 2, 0, 1 ]$&

$
\begin{array}{l}
      y^2 + \omega^5 y + 1 \\
      y^2 + \omega^7 y + 1 \\

      x y + \omega x + \omega^3 y + 1 \\
      x y + \omega^3 x + \omega y + 1 \\
      x^2 + \omega^5 x + 1 \\
      x^2 + \omega^7 x + 1 \\
\end{array}$
&
$[ 2, 1, 0, 2, 2, 1 ]$&

$
\begin{array}{l}
      y + \omega,  y + \omega^3 \\
      y + \omega^5, y + \omega^7 \\
      x + \omega, x + \omega^3 \\
      x + \omega^5,x + \omega^7 \\

      x y + \omega^5 x + \omega^7 y + 1 \\
      x y + \omega^7 x + \omega^5 y + 1 \\
\end{array}$\\
\hline
\end{tabular}
\end{table}

Finally, we deal with fractional functions of type $\frac{\widetilde{h}}{h}$, for odd characteristic $p$ and $\deg(h)=3$. In particular, we are able to provide 
 conditions on the coefficients of $h$ which ensures that  $\frac{\widetilde{h}}{h}$ is a permutation of $\mu_{q+1}$. This could be useful to investigate the permutation property of the associated polynomials.
\begin{proposition}\label{Finale} Let $p>2$.
Let $A,B,C \in \mathbb{F}_q$, $A,B\neq 0$. Let $h_{A,B,C}(x)={B x^3+Cx^2+x+A}$ The function $f_{A,B,C}(x) = \frac{\widetilde{h}_{A,B,C}(x)}{{h}_{A,B,C}(x)}$ permutes $\mu_{q+1}$ in the following cases:
\begin{itemize}
\item  $A^2 - AC - B^2 + B=0$, and $AT^2 + (1-B)T +A\in \mathbb F_q[T]$ has distinct roots in $\mathbb{F}_{q}$;
\item $q \equiv 1 \pmod 3$,  $B=(3 A C  + C^2 - 1)/3$, and the polynomial $3T^2    -3(3A+C)T + (3A+C)^2   - 1\in \mathbb F_q[T]$ has distinct roots in $\mathbb{F}_q$. 
\end{itemize}
\end{proposition}
\proof
\begin{itemize}
\item If $p$ is odd and $A^2 - AC - B^2 + B=0$ then the equation of the associated  curve $\mathcal{C}_{f_{A,B,C}}$ can be written as  
$$(Ax^2 + (1-B)x +A)(Ay^2 + (1-B)y +A)=0.$$
Let $k_1$ and $k_2$ the distinct roots of  $AT^2 + (1-B)T +A$; then the curve splits into the four lines  $x=k_1$, $x=k_2$, $y=k_1$, and $y=k_2$. Since $k_1,k_2 \in \mathbb F_q$, these lines do not have points in $\mu_{q+1}^2$. In fact, $\mathbb{F}_q\cap \mu_{q+1}=\{\pm 1\}$ and $k_1,k_2=\pm 1$ would imply $k_1=k_2$.
%
\item Let $p\neq 3$ and assume that $B=(3 A C  + C^2 - 1)/3$. Then the equation of the associated curve $\mathcal{C}_{f_{A,B,C}}$ can be written as
$$(x y+\alpha x+(3 A+C-\alpha)y+1)(x y+(3 A+C-\alpha) x+\alpha y+1)=0$$
where $\alpha$ satisfies 
$$3\alpha^2    -3(3A+C)\alpha + (3A+C)^2   - 1=0.$$
By assumption $\alpha \in \mathbb{F}_q$. 
The conic $x y+\alpha x+(3 A+C-\alpha)y+1=0$ has a point $(x_0,y_0)$ in $\mu_{q+1}^2$ if and only if 
$$y_0^{q}= \left(-\frac{\alpha x_0+1}{x_0+3 A+C-\alpha}\right)^q
=-\frac{\alpha+x_0}{1+(3 A+C-\alpha)x_0}$$
equals $1/y_0=\frac{x_0+3 A+C-\alpha}{\alpha x_0+1}$. Then 
$$(\alpha+x_0)(\alpha x_0+1)-((3 A+C-\alpha)+x_0) ((3 A+C-\alpha) x_0+1)=0$$
and hence
$$(3A+C-2\alpha)(x_0^2+(3A+C)x_0 + 1)=0$$
for some $x_0\in \mu_{q+1}$.
Since the roots of $F(T)=3T^2    -3(3A+C)T + (3A+C)^2   - 1$ are distinct, $3A+C-2\alpha=0$ cannot hold. Also, the roots of $T^2+(3A+C)T + 1=0$ are in $\mathbb{F}_q$ since $q\equiv 1 \pmod 3$ ensures that $-3$ is a square in $\mathbb F_q$ (see e.g. \cite[Lemma 4.5]{giul}); hence, $x_0\in \{\pm 1\}$. This implies $3A+C=\pm 2$ and $F(T)=3(T\pm1)^2$, a contradiction.

\end{itemize}
\endproof

\begin{corollary}\label{Cor:Finale}
Let $p>2$, $A,B,C \in \mathbb{F}_q$, $A,B\neq 0$. The polynomial $F_{A,B,C}(x) = Ax^{q(q-1)+1}+Bx^{2(q-1)+1}+Cx^{q-1+1}+x$ permutes $\mathbb{F}_{q^2}$ if either 
\begin{itemize}
\item[(i)]  $A^2 - AC - B^2 + B=0$ and $AT^2 + (1-B)T +A$ has distinct roots in $\mathbb{F}_{q}$,
\end{itemize}
or
\begin{itemize}
\item[(ii)] $q \equiv 1 \pmod 3$,  $B=(3 A C  + C^2 - 1)/3$, and the polynomial $F(T)=3T^2    -3(3A+C)T + (3A+C)^2   - 1$ has distinct roots in $\mathbb{F}_q$. 
\end{itemize}

\end{corollary}
\proof
By Theorem \ref{ThPLZ} $F_{A,B,C}(x)$ permutes $\mathbb{F}_{q^2}$ if and only if $x(Ax^{q}+Bx^{2}+Cx+1)^{q-1}$ permutes $\mu_{q+1}$. On $\mu_{q+1}$ the previous function equals $\frac{A x^3+x^2+C x+B}{B x^3+Cx^2+x+A}$. By Proposition \ref{Finale} the assertion follows.
\endproof

We now show that the second condition of Corollary \ref{Cor:Finale}$(i)$ can be replaced by
$\frac{B-1}{A}\notin \mathbb{F}_q$.

\begin{proposition}
Let $A,B,C \in \mathbb{F}_{q^2}$, $A,B\neq 0$. If $A^2 - AC - B^2 + B=0$ and $\frac{B-1}{A}\notin \mathbb{F}_q$ then the polynomial $F_{A,B,C}(x) = Ax^{q(q-1)+1}+Bx^{2(q-1)+1}+Cx^{q-1+1}+x$ permutes $\mathbb{F}_{q^2}$.

\end{proposition}
\proof
Let $f_{A,B,C}$ be the fractional polynomial associated with the polynomial $F_{A,B,C}$. We know by Proposition \ref{Finale} that in this case the equation of the corresponding curve can be written as $(Ax^2 + (1-B)x +A)(Ay^2 + (1-B)y +A)=0$.  Note that the roots of $AT^2 + (1-B)T +A$ are of type $\{z,1/z\}$ and therefore they belong to  $\mu_{q+1}$ if and only if $(B-1)/A=-z-z^{-1}=-z-z^q$ for some $z\in \mu_{q+1}$. Note that $z+z^q\in \mathbb{F}_{q}$ and by assumption this is not possible. 
\endproof


\section{Acknowledgment}
The authors were supported in part by Ministry for Education, University and Research of Italy (MIUR) (Project PRIN 2012 Geometrie di Galois e strutture di incidenza") and by the Italian National Group for Algebraic and Geometric Structures and their Applications (GNSAGA - INdAM).

\end{document}